\documentclass[letterpaper, 10 pt, conference]{ieeeconf}
\IEEEoverridecommandlockouts  
\usepackage{lineno,setspace, nicefrac, subcaption, caption}
\usepackage{mathrsfs}
\usepackage{tikz}
\usetikzlibrary{automata,arrows,positioning,calc}
\usepackage[mathscr]{eucal}

\DeclareGraphicsExtensions{.png,.jpeg,.jpg,.pdf}
\graphicspath{{./Figures/}}

\newtheorem{theorem}{Theorem}[section]

\newtheorem{proposition}[theorem]{Proposition}

\newtheorem{lemma}[theorem]{Lemma}

\newtheorem{remark}[theorem]{Remark}

\newcommand{\lb}{\left\{}
\newcommand{\rb}{\right\}}
\newcommand{\set}[1]{\left\{ #1 \right\}}

\newcommand{\Def}{\overset{\textbf{def}}{=}}
\newcommand{\RR}{\mathbb{R}}

\newcommand{\eps}{\varepsilon}
\newcommand{\bOne}{\textcolor{black}{\mathbf{1}}}

\newcommand{\MID}{\mathcal{T}}

\newcommand{\abs}[1]{\lvert #1 \rvert}
\newcommand{\norm}[1]{\left\| #1 \right\|}

\newcommand{\sT}{\mathsf T}

\newcommand{\bigmid}{\bigl\vert}

\usepackage{amsopn}

\DeclareMathOperator{\ball}{Ball}

\usepackage{mathptmx} 
\usepackage{times} 
\usepackage{amsmath} 
\usepackage{amssymb}  

\usepackage{color}
\definecolor{darkgreen}{rgb}{0,0.5,0}
\usepackage[colorinlistoftodos]{todonotes}
\usepackage{soul}
\setulcolor{red} 
\setstcolor{red} 
\sethlcolor{yellow} 

\usepackage{dirtytalk}

\title{\LARGE \bf Near Collision and Controllability Analysis of Nonlinear Optimal Velocity Follow-the-Leader Dynamical Model In Traffic Flow}
\author{Hossein Nick Zinat Matin$^{1}$ and {Maria Laura} {Delle Monache} $^{1}$
\thanks{$^{1}$ Hossein Nick Zinat Matin and Maria Laura Delle Monache are with the Department of Civil and Environmental Engineering,
        University of California, Berkeley,
        {\tt\small h-matin@berkeley.edu, mldellemonache@berkeley.edu}}
}

\date{\today.}

\begin{document}

\maketitle 
\thispagestyle{empty}
\pagestyle{empty}

\begin{abstract}
This paper examines the optimal velocity follow-the-leader dynamics, a microscopic traffic model, and explores different aspects of the dynamical model, with particular emphasis on collision analysis. More precisely, we present a rigorous boundary-layer analysis of the model which provides a careful understanding of the behavior of the dynamics in trade-off with the singularity of the model at collision which is essential in the controllability of the system. 

\end{abstract}
\section{Introduction and Related Works}
The emergence of autonomous driving technologies such as adaptive cruise control and self-driving systems has created different theoretical challenges in modeling and analysis of the governing dynamics of the traffic flow.

Traffic flow dynamics has been a widely studied research area for decades, with literature devoted to various models based on macroscopic, mesoscopic, and microscopic descriptions of traffic flow \cite{treiber2013traffic}. The microscopic class of dynamics considers individual vehicles and their interaction. The earliest car following models date back to the works of \cite{pipes1953operational,newell1961nonlinear, gazis1959car, chandler1958traffic}. Nonlinear follow-the-leader dynamics can be traced back to \cite{herman1959single} and \cite{gazis1961nonlinear} among others. The celebrated Optimal Velocity (OV) dynamical model was introduced and analyzed in \cite{bando1995dynamical, bando1994structure, bando1998analysis} and numerous following studies. In this paper, we consider the Optimal Velocity Follow-the-Leader (OVFL) dynamical model which is shown to possess favorable properties both from a practical and theoretical point of view \cite{tordeux2014collision, stern2018dissipation,nick2022near,matin2020nonlinear,matin2019ConvergenceRate,dellemonache2019pardalos}. 

The optimal velocity part of the OVFL model with a positive coefficient defines a target velocity based on the distance between each vehicle and its preceding one. Comparing the target and the current velocities, the acceleration/deceleration will be encouraged by the OV model. The follow-the-leader term explains the force that tries to match the vehicle's velocity with the preceding one. As the instantaneous relaxation time (i.e. $(x_{n-1} - x_n)/\beta$ in \eqref{E:initial_model}) decreases, a singularity occurs at collision. Understanding the interaction between such a singularity and the behavior of OVFL dynamics near collision is the main focus of this paper.

Stability analysis of platoon of vehicles following \eqref{E:initial_model} has been studied from various points of view such as string stability, \cite{chandler1958traffic, bando1998analysis, kometani1958stability, wilson2011car, gunter2020commercially, giammarino2021traffic, piu2022stability}. Analysis of collision has been addressed from different standing points in some prior works. In a simulation-based study, \cite{davis2003modifications} investigates the likelihood of collision as a consequence of drivers' reaction time. A Lyapunov-based analysis in a neighborhood of the equilibrium point has been studied in \cite{davis2014nonlinear, tordeux2014collision}. Nonlinear stability analysis and collision avoidance based on the safe distance is studied in \cite{magnetti2021nonlinear} for the OV model. 

\textbf{Focus and Contribution.} In contrast to the stability-based analysis of OVFL dynamics, in this paper, we are interested in the analysis of collision (e.g. in a platoon of connected autonomous vehicles which are governed by such a dynamical model). In other words, our main focus is on understanding the interplay between the behavior of the OVFL dynamics and the singularity introduced in \eqref{E:initial_model} at collision, through a careful and mathematically rigorous investigation. 

Our boundary-layer analysis results are strongly dependent on the initial values which allow us to study the effect of singularity when the vehicles are in a near-collision region. Such analytical understanding is crucial in analyzing the behavior of the system in real-world conditions such as in the presence of noise and perturbation. In such conditions, sooner or later any physical system will be pushed into various states. Therefore it is necessary and insightful to understand the deterministic behavior of the system in the proximity of critical states. 

As a consequence of our analysis, we show that the collision in the system does not happen and hence the system is well-posed. In addition, our analysis applies to multiple-vehicle which extends the results of \cite{nick2022near}. 

The organization of this paper is as follows. We start by introducing the dynamical model. Then, we prove some essential properties of the dynamics between the first two vehicles which will be used in the analysis of the other following vehicles. Then, we study the behavior of the trajectory of other vehicles with respect to that of the first two and we show the main result of the paper.
\section{Mathematical Model}
\begin{figure}
    \centering
    \includegraphics[width=3.3in]{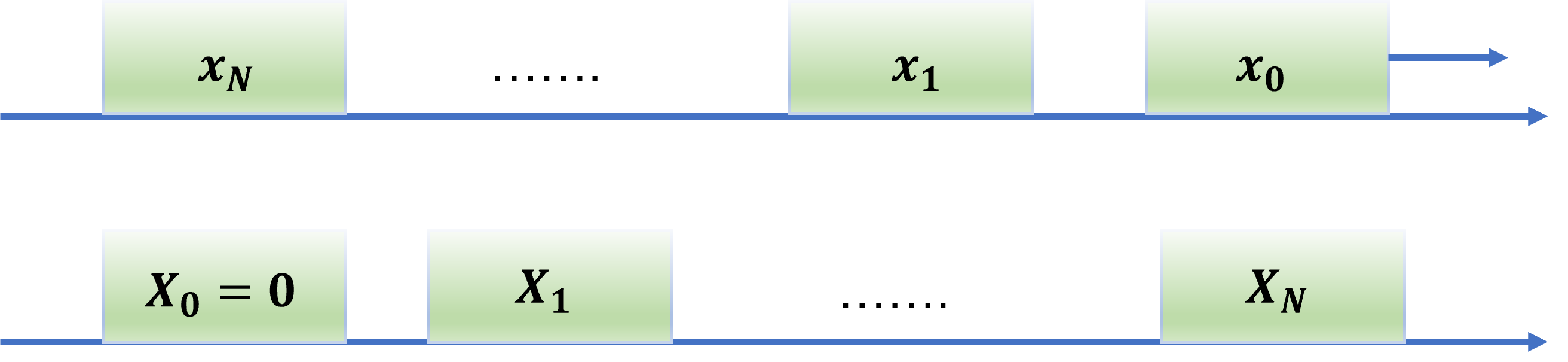}
    \caption{\small{The first illustration shows the position and direction of the vehicles. The second illustration depicts the relative position $X_n = x_0 - x_n$.}}
    \label{fig:direction}
\end{figure}

\begin{figure}
    \centering
    \includegraphics[width=3in]{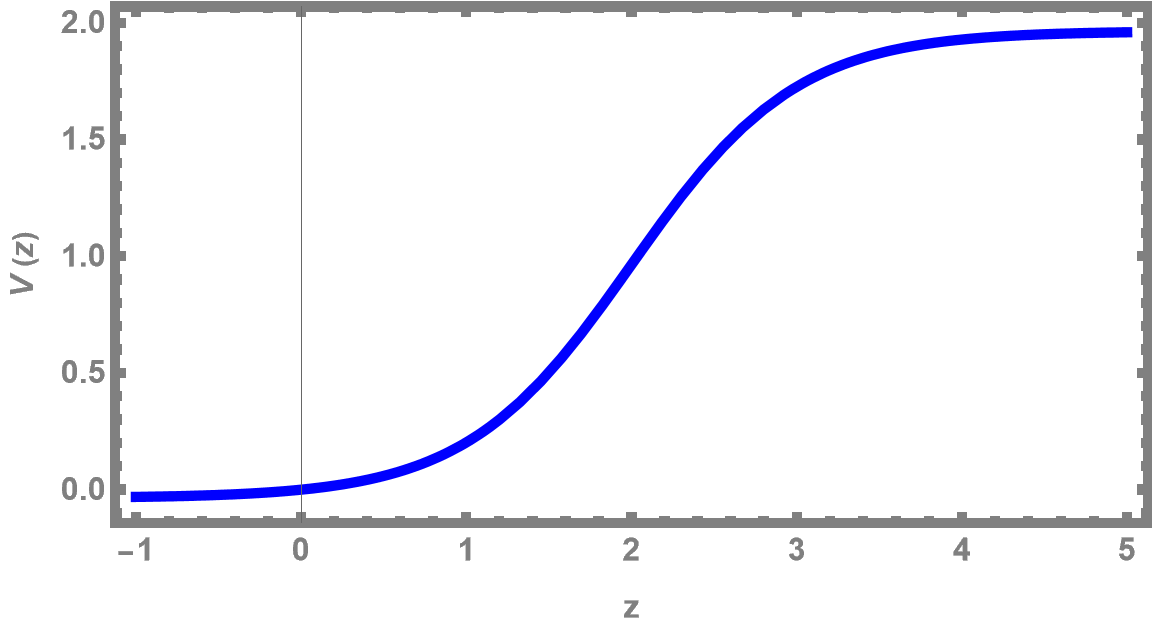}
    \caption{Function $V$ in \eqref{E:V_function}.}
    \label{fig:Vplot}
\end{figure}

We consider $N+1$ number of vehicles and each vehicle $n = 0,1, \cdots, N$ has position $x_n$ and velocity $y_n = \dot x_n$ such that $x_N < x_{N-1} < \cdots< x_0$ (see Figure \ref{fig:direction}).
We assume that the first vehicle is moving with a constant velocity $\bar v$; i.e. $\dot x_0(t) = y_0(t) = \bar v$ with the initial value $(x_\circ, y_\circ)$. The OVFL model for $n \ge 1$ can be presented in the form of
\begin{equation}
    \label{E:initial_model}
    \begin{cases}
    \dot x_n(t) = y_n(t) \\
    \dot y_n(t) = \alpha \lb V(x_{n -1} - x_n) + y_n(t)) \rb + \beta \frac{y_{n -1} - y_n }{(x_{n -1} - x_n)^2} \\
    (x_n(0), y_n(0)) = (x_{n, \circ}, y_{n, \circ})
\end{cases}\end{equation}
where function $V$ is monotonically increasing, bounded, and Lipschitz continuous function. 
In this paper, we consider 
\begin{equation}
    \label{E:V_function}
    V(x) \Def \tanh(x - 2) - \tanh(-2);
\end{equation}
as illustrated in Figure \ref{fig:Vplot}. 
We define $V(\infty) = 1.96$  as the scaled maximum possible speed. It should also be noted that $\Delta x_n = x_n - x_{n-1} = l$, where $l$ is the length of vehicles, is interpreted as the collision between vehicle $n-1$ and $n$. In this paper, without loss of any generality, we drop the constant $l$ and consider $\Delta x_n =  0$ as a collision. 
Therefore, 
\begin{equation}\label{E:V_infty}
y_n(t) \in [0, V(\infty)), \quad  n \in \set{0, \cdots, N}, \quad t \ge 0. \end{equation}
For simplicity of the following analysis and interpretation of the results, we define a Galilean change of variable in \eqref{E:initial_model}
\begin{equation*}\begin{split}
    X_n(t) \Def x_0(t) - x_n(t), \quad n =1, \cdots, N \\
    Y_n(t) \Def y_0(t) - y_n(t) = \bar v - y_n(t), \quad n = 1 , \cdots, N
\end{split}\end{equation*}
Consequently, the dynamics of \eqref{E:initial_model} can be rewritten as
\begin{equation}
    \label{E:main_dynamics}
    \begin{cases}
    \dot X_n(t) = Y_n(t) \\
    \dot Y_n(t) = -\alpha \lb V( X_n- X_{n -1}) + Y_n(t) - \bar v \rb -\beta \frac{Y_n -Y_{n -1}}{ (X_n-X_{n -1})^2} \\
    (X_n(0), Y_n(0)) = (X_{n, \circ}, Y_{n, \circ})
\end{cases}\end{equation}
for $n \in \set{1, \cdots, N}$ with the convention that $X_0 = Y_0 = 0$. It should be noted that in \eqref{E:main_dynamics} we have that $X_N > \cdots> X_1 > X_0 = 0$; see Figure \ref{fig:direction}. 
In addition, following \eqref{E:V_infty}, we have that 
\begin{equation}
    \label{E:speed_limit}
    Y_n(t) \in (\bar v - V(\infty), \bar v], \quad t \ge 0.
\end{equation}
This is in particular important in choosing the initial values of the dynamics. 
The dynamical system \eqref{E:main_dynamics} has a unique equilibrium solution when all the vehicles are equidistantly located and moving with the same velocity \cite{bando1998analysis, piu2022stability}. Mathematically, for each $n \in \set{1, \cdots, N}$
\begin{equation*}
    (X_n^\infty , Y_n^\infty) = (n V^{-1}(\bar v) , 0)= (n X_\infty, 0).
\end{equation*}
where 
\begin{equation*}
        X_\infty =  V^{-1}(\bar v)=2 + \tanh^{-1}(v_\circ + \tanh(-2)),
\end{equation*}

\section{Dynamics of the First Two Vehicles}\label{S:two_vechicles}
The behavior of the dynamics propagates from the leading vehicles to the following ones. Therefore, we need to start by understanding the interaction between the first two vehicles. 
\subsection{Hamiltonian and Boundedness of Solution} 
\begin{figure}
    \centering
    \includegraphics[width=0.8\columnwidth]{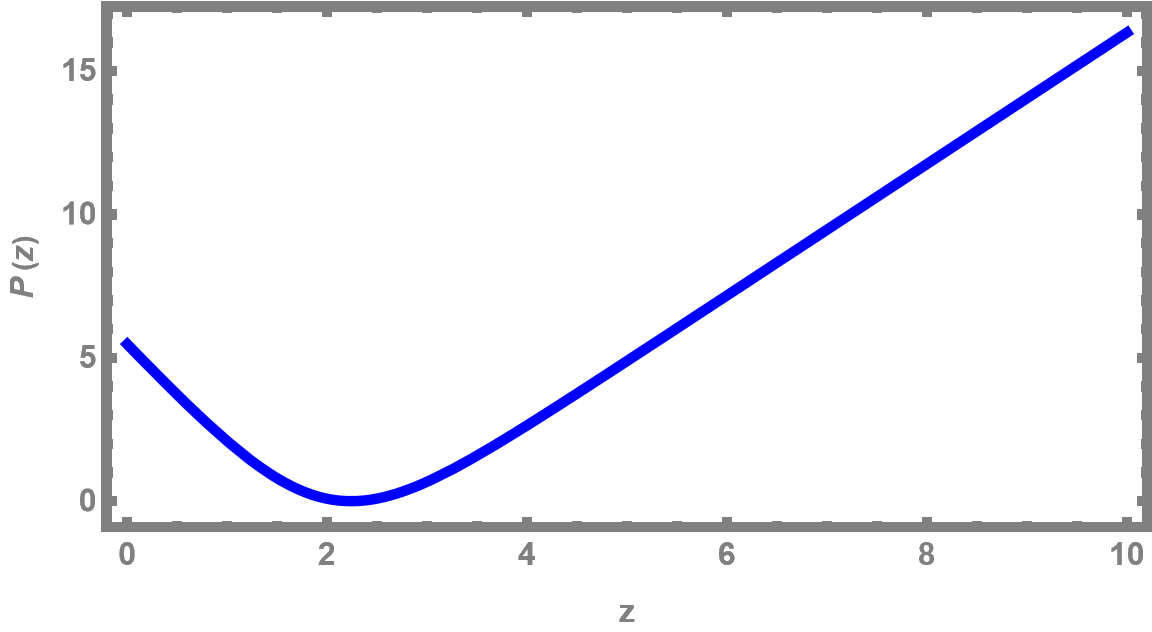}
    \includegraphics[width=0.8\columnwidth]{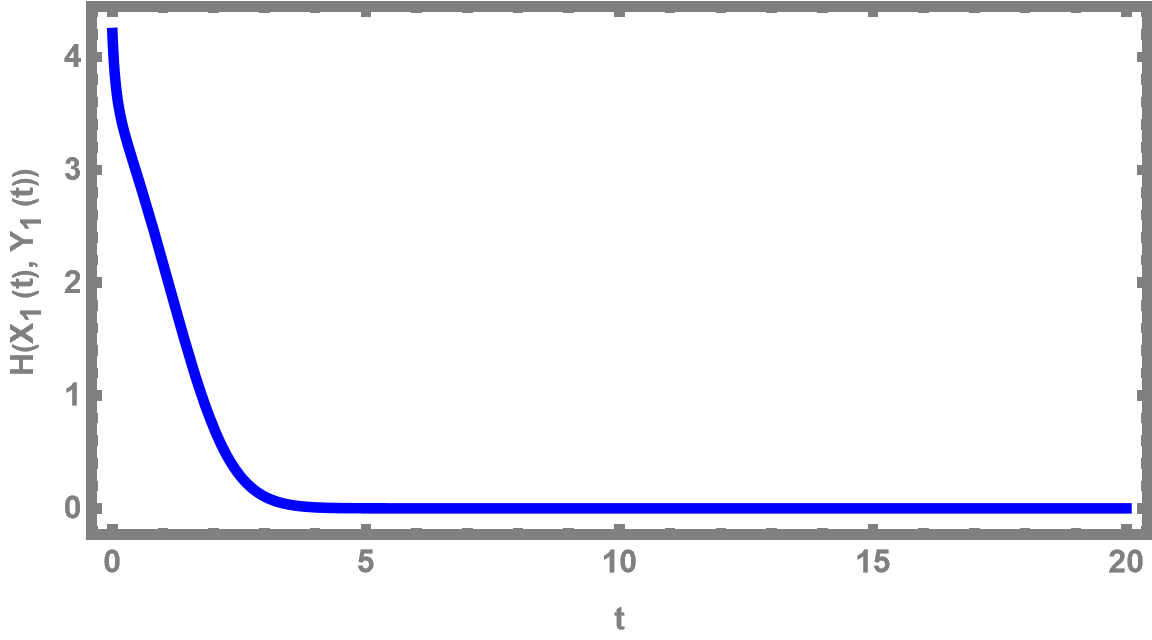}
    \caption{Potential $P$ \eqref{eq:potential} and Hamiltonian $H$ \eqref{E:Hamiltonian} along the trajectory $(X_1(t), Y_1(t))$.}
    \label{fig:geometry}
\end{figure}

\begin{figure}
    \centering
    \includegraphics[width=.8\columnwidth]{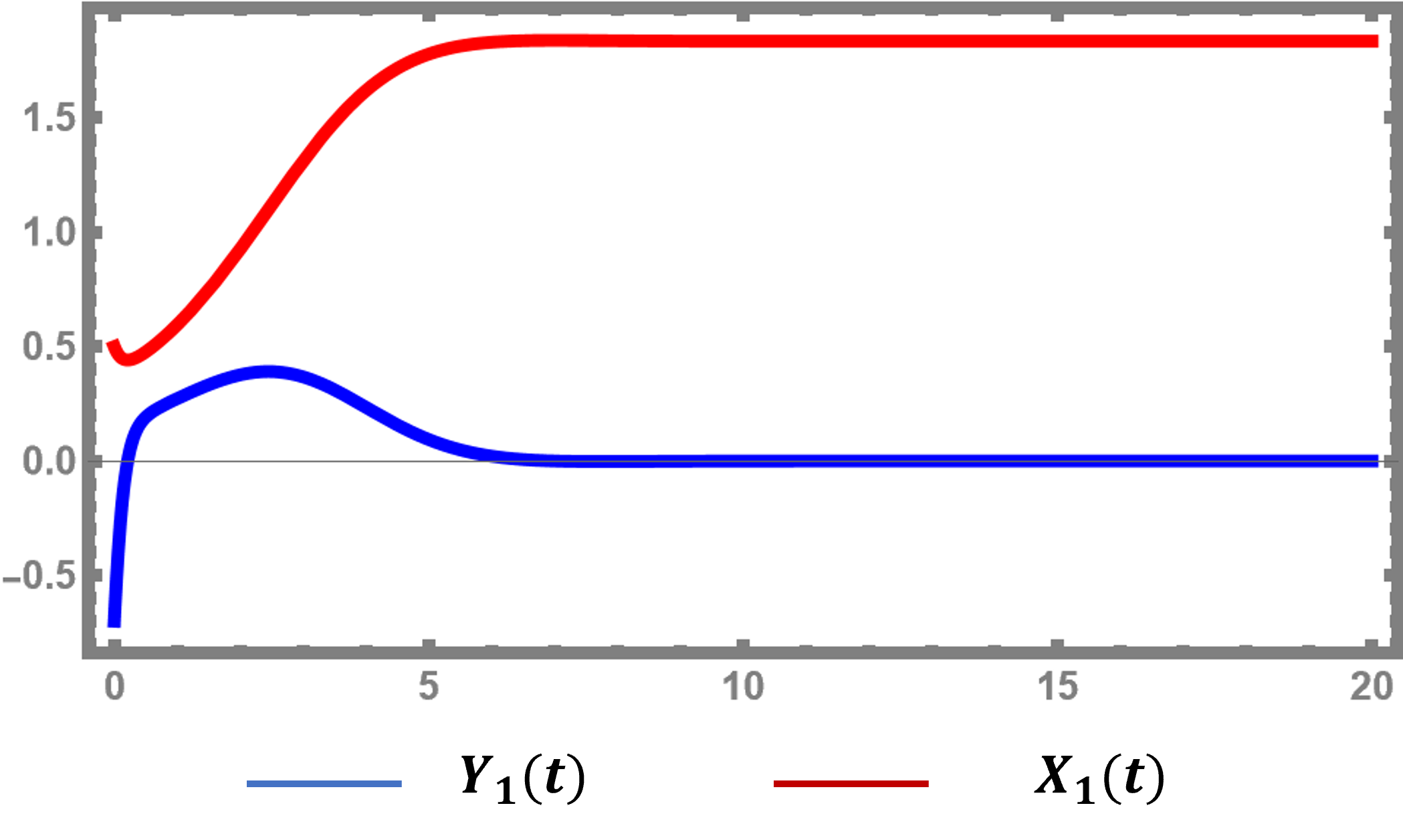}
    \includegraphics[width=.8\columnwidth]{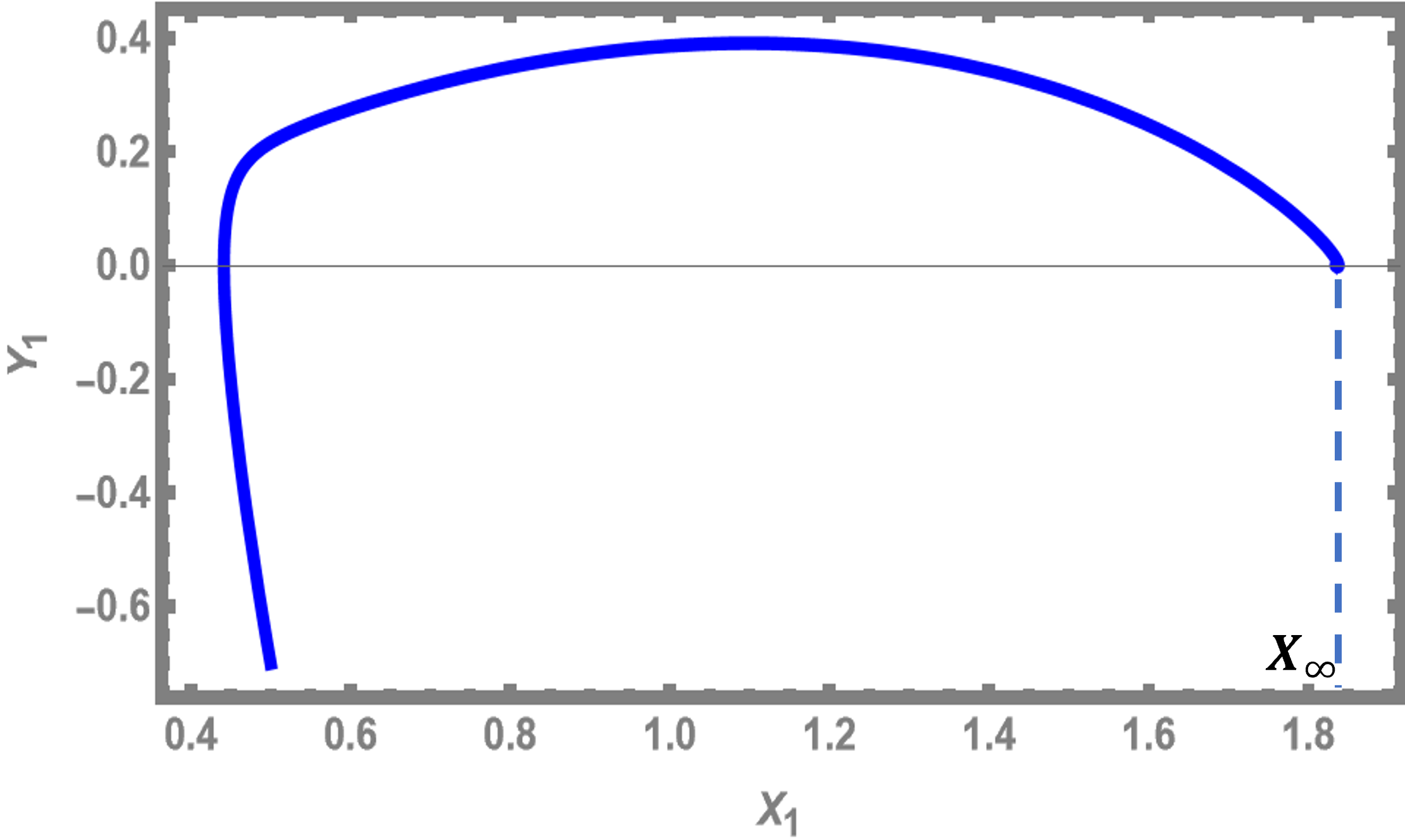}
    \caption{\small{Trajectory of the first two vehicles. The top plot shows the trajectory of $t \mapsto X_1(t)$ and $ t \mapsto Y_1(t))$ separately. The bottom plot show the orbit of this dynamic for $\alpha = 2$, $\beta = 1$, $\bar v = 0.8$, $(X_{1, \circ}, Y_{1, \circ}) = (0.5, -0.7)$.}}
    \label{fig:4-5}
\end{figure}

In this section, we consider $N =1$ in  \eqref{E:main_dynamics}, the dynamics between the first two vehicles. Following \cite{nick2022near}, first we recall few properties of the dynamics of $(X_1(t), Y_1(t))$. The main properties of the dynamics of $(X_1(t), Y_1(t))$ can be obtained by defining the Hamiltonian function
\begin{equation}\label{E:Hamiltonian}
    H(x, y) \Def \frac 12 y^2 + P(x) \qquad (x, y) \in \RR^2
\end{equation}
with the potential function
\begin{equation}\label{eq:potential}
    P(x) \Def \alpha \int_{x'=X_\infty }^x \lb V \left(x' \right) - \bar v \rb dx'; \qquad x \in \RR.
\end{equation}
The graphical depiction of \eqref{E:Hamiltonian} and \eqref{eq:potential} are shown in Figure \ref{fig:geometry}. Since 
\begin{equation*}
    P'(x) = \alpha \lb V (x) - \bar v \rb, \qquad x>0
\end{equation*}
we can write the dynamics \eqref{E:main_dynamics} for $n =1$ as a \emph{damped Hamiltonian system}
\begin{equation}\label{E:main_deterministic} \begin{aligned} \dot X_1(t)&=Y_1(t)=\frac{\partial H}{\partial y}(X_1(t), Y_1(t)) \\
\dot Y_1(t) &= -P'(X_1(t))-\alpha Y_1(t) -\beta\frac{Y_1(t)}{X_1^2(t)} \\
&= -\frac{\partial H}{\partial x}(X_1(t), Y_1(t))-\lb \alpha +  \frac{\beta}{X_1^2(t)}\rb Y_1(t) \end{aligned} \end{equation}
for $t\ge 0$.
\begin{remark}\label{R:Pstructure} The function $V$ is strictly increasing which implies that $V(x)-v_\circ <0$ if $x < X_\infty$ and $V(x)-v_\circ > 0$ if $x>X_\infty$. Thus $P$ will be increasing on the interval $(X_\infty,\infty)$, decreasing on the interval $(0,X_\infty)$ and its minimum point is located at $X_\infty$.  Furthermore $\lim_{x\nearrow \infty}P(x)=\infty$. \end{remark}
Using $H(x, y)$, we can show that the solution $(X_1, Y_1)$ is bounded. In particular, we have that
\begin{equation}\label{E:energyLevel}
\begin{split}
    \dot H(X_1(t), Y_1(t))&= \dot X_1(t) \frac{\partial H}{\partial x}(X_1(t), Y_1(t)) \\
    & \qquad \qquad + \dot Y_1(t) \frac{\partial H}{\partial y}(X_1(t), Y_1(t))\\
    &= -\lb \alpha + \frac{\beta}{X_1(t)^2}\rb Y_1(t)^2 \le 0
\end{split}\end{equation}
See the illustration of function $H(X_1(t), Y_1(t))$ in Figure \ref{fig:geometry}. This can be interpreted as decreasing energy in the system. Now, if we set
\begin{equation}\label{E:hcircdef} h_\circ \Def H(X_{1,\circ},Y_{1, \circ}), \end{equation}
then from the definition of $H$, we have that 
\begin{equation}
    \tfrac 12 Y_1^2(t) \le H(X_1(t), Y_1(t)) \le h_\circ
\end{equation}
where the second inequality follows from \eqref{E:energyLevel} and \eqref{E:hcircdef}. Therefore, 
\begin{equation}\label{E:Y_bound}
    \abs{Y_1(t)} \le \bar y\Def \sqrt{2 h_\circ}, \quad t \ge 0.
\end{equation}
Similarly, we can show by increasing behavior of $P$ on $(X_\infty, \infty)$ that 
\begin{equation}\label{E:X1_bound}
    X_1(t) \le \bar x \Def h_\circ, \quad t \ge 0.
\end{equation}
Finally, it is shown that the solution $X_1(t)$ never hit zero which can be interpreted as no collision between the first two vehicles.
In particular, 
\begin{equation}\label{E:X1_LB}
X_1(t) \ge \delta_1, \quad t \ge 0
\end{equation}
where, the lowerbound $\delta_1 = \delta_1(X_{1, \circ}, Y_{1, \circ})$ depends on the initial values of the system. Figure \ref{fig:4-5} shows the trajectory of the dynamical model for the first two vehicles.

Now, we need to develop some more properties of the trajectory of flow $(X_1(t), Y_1(t))$. 
\begin{remark}\label{R:global_asym_stab}
    Using \eqref{E:Y_bound} and \eqref{E:X1_LB}, set $(0, \bar x] \times [-\bar y, \bar y]$ is positively invariant with respect to the flow $(X_1(t), Y_1(t))$ and has a compact closure. Furthermore, \eqref{E:energyLevel} precludes a periodic orbit. Therefore, an application of the Poincar\'e-Bendixon theorem \cite[Section 7.3]{teschl2012ordinary}, suggests that the equilibrium solution $(X_\infty, 0)$ is globally asymptotically stable.
\end{remark}

\begin{remark} [\textbf{Initial Data}] In this paper, we are mainly concerned with the boundary-layer analysis of the system near collision (the situation that can happen, for instance, as a result of instantaneous perturbation in the system; like sudden braking of the leading vehicles which propagates). In other words, we are interested in the case that the distance between the corresponding consecutive vehicles becomes relatively small. In particular, we consider the initial values $\Delta X_{n, \circ} < X_\infty$, for the respective $n \in \set{1, \cdots, N}$. 

In addition, suppose that $Y_{1, \circ} <0$. Fix a time $T>0$. Since $X_{1, \circ} <X_\infty$, the dynamics of \eqref{E:main_dynamics} for $N =1$ suggest that $\dot Y_1(t) >0$ for $t\in (0, \eps_\circ)$; some neighborhood of time zero. On the other hand, since the $X_{1, \circ}$ is relatively small, the dominant term in the dynamics of $\dot Y_1$ in \eqref{E:main_dynamics} is $-\nicefrac{\beta Y_1(t)}{(X_1(t))^2}$ for $t \in (0, \eps_\circ)$. Hence, for sufficiently large $\beta$, $Y_1(\bar t) >0$ for some $\bar t <T$ (see Figure \ref{fig:4-5}). Therefore, in this paper, without sacrificing any generality, it is sufficient to consider $Y_{1, \circ} >0$; otherwise, the same analysis follows after shifting the initial time to $\bar t$. 

Moreover, this assumption will not affect the generality of our follower vehicles' analysis in the next section. More precisely, let $N = 2$. The interaction between two consecutive vehicles depends on their relative speed, i.e. $Y_2 - Y_1$ (rather than merely the relative velocity $Y_1$ of the leading vehicle) which will be analyzed in its full generality. In particular, as we will see, the most interesting case for the purpose of our boundary layer analysis will be $Y_n - Y_{n-1} <0$, $n \ge 2$, which implies that the following vehicle is moving faster than the leading one. This can potentially result in a collision. We will discuss this case in detail in the next section. 

\end{remark} 
\subsection{Controlling the Behavior of the Dynamics by Controlling the Parameters}
In this section, we study the behavior of the trajectory of $t \mapsto Y_1(t)$ for $(X_{1, \circ}, Y_{1, \circ})\in (0, X_\infty] \times \RR_+$. We define
\begin{equation}
    \label{E:equil_time}
    \MID_\infty \Def \inf \set{t \ge 0: X_1(t) = X_\infty}
\end{equation}
as the first time for which the trajectory $t \mapsto X_1(t)$, starting from the initial data $(X_{1, \circ}, Y_{1, \circ})$, approaches $X_\infty$. The change of variables $u \Def x - X_\infty$ and $v \Def y$ help us standardize the stability analysis by translating the equilibrium point to the origin. 
The Hamiltonian can be rewritten as 
\begin{equation}\label{E:new_hamiltonian}\begin{split}
        H(u , v) &= \tfrac 12 v^2 + \tilde P(u) \\
        & \Def \tfrac 12 v^2 + \alpha \int_0^u \lb V(u' + X_\infty) - \bar v \rb du' \\
        \frac{dH}{dt}(u(t), v(t)) &= - \left(\alpha + \frac{\beta}{(u(t) + X_\infty)^2} \right) v^2(t).
    \end{split}\end{equation}
The main result of this section expresses that by controlling the parameters $\alpha$ and $\beta$ we can control the behavior of the trajectory $t \mapsto Y_1(t)$. In particular,
\begin{theorem}
        \label{T:exp_stable_convergence}
        Starting from $(X_{1, \circ}, Y_{1, \circ})\in (0, X_\infty] \times \RR_+$, for sufficiently large values of $\alpha$ and $\beta$, we have that 
        \begin{equation*}
            \limsup_{t \nearrow \MID_\infty } v(t) = 0.
        \end{equation*}
    \end{theorem}
    This implies that by controlling $\alpha$ and $\beta$, the flow $Y_1(t)$ will be absorbed to the equilibrium point as $t \nearrow \MID_\infty$. We postpone the proof of this theorem until some preliminary results are established. The following  lemmas explain the behavior of the trajectory $t \mapsto Y_1(t)$ around the equilibrium point. Figure \ref{fig:9-10} is provided as a graphical aide to the proofs. 
\begin{lemma}\label{T:stay_positive}
The set $U_1 \Def \set{ (x, y): x < X_\infty, y\in \RR_+}$ is invariant with respect to the trajectory $[0, \MID_\infty) \ni t \mapsto (X_1(t), Y_1(t))$. In other words, if $(X_{1,\circ}, Y_{1, \circ}) \in U_1$, then $Y_1(t) >0$ for $t \in [0, \MID_\infty)$.
\end{lemma}
\begin{proof}
We use the proof by contradiction to show the result. Suppose that there exists a time $t_\circ < \MID_\infty$ such that $Y_1(t_\circ) = 0$. By the continuity of the solution, we have that $\dot Y_1(t_\circ) <0$. On the other hand, 
    \begin{equation*}
        \begin{split}
            \dot Y_1(t_\circ) &= - \alpha \lb V(X_1(t_\circ)) - \bar v + Y_1(t_\circ) \rb - \beta \frac{Y_1(t_\circ)}{(X_1(t_\circ))^2}\\
            & = -\alpha \lb V(X_1(t_\circ)) - \bar v \rb  >0
        \end{split}
    \end{equation*}
    where the last inequality holds since $X_1(t_\circ) < X_\infty$. But this is a contradiction and hence the result follows (see Figure \ref{fig:9-10}). 
    \begin{figure}
        \centering
        \includegraphics[width= 3.5in]{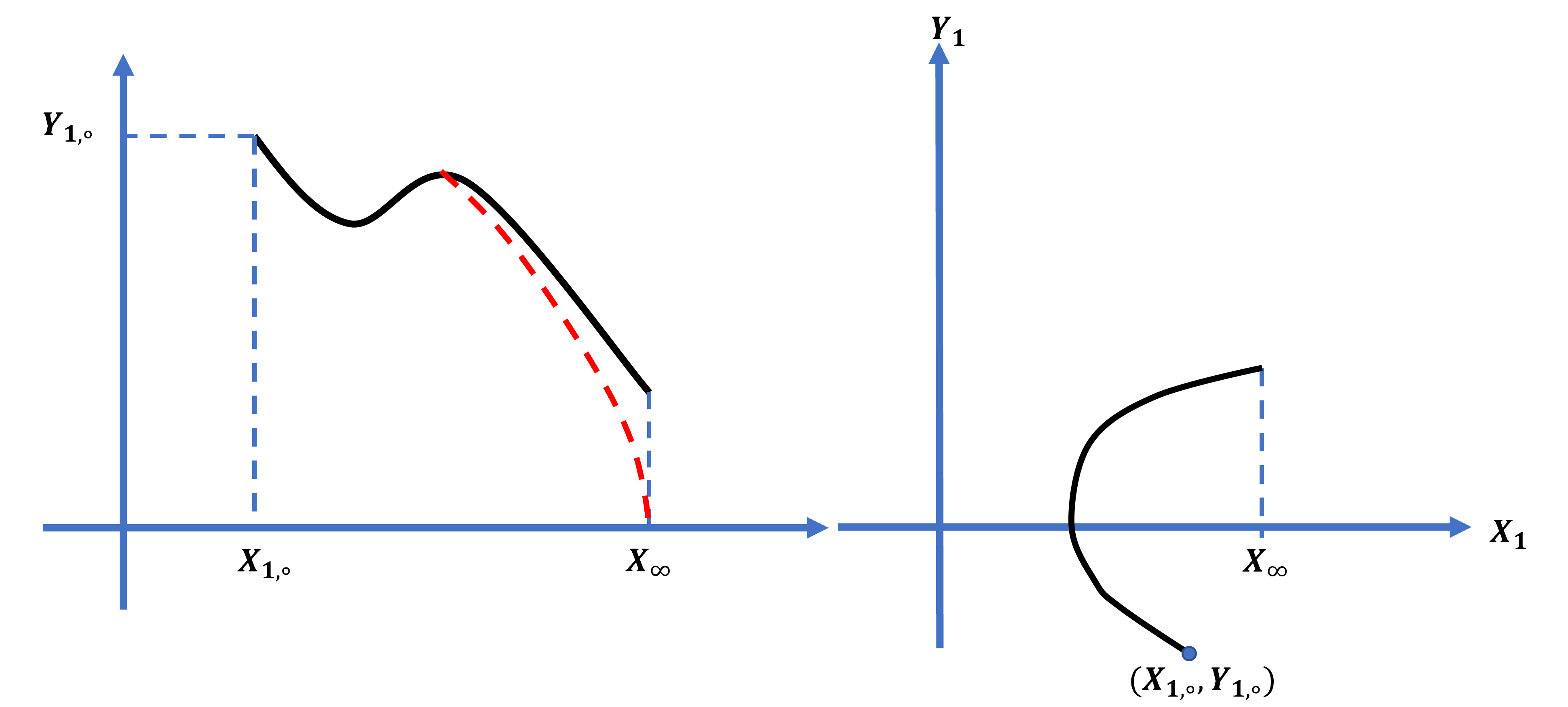}
        \includegraphics[width=3.5in]{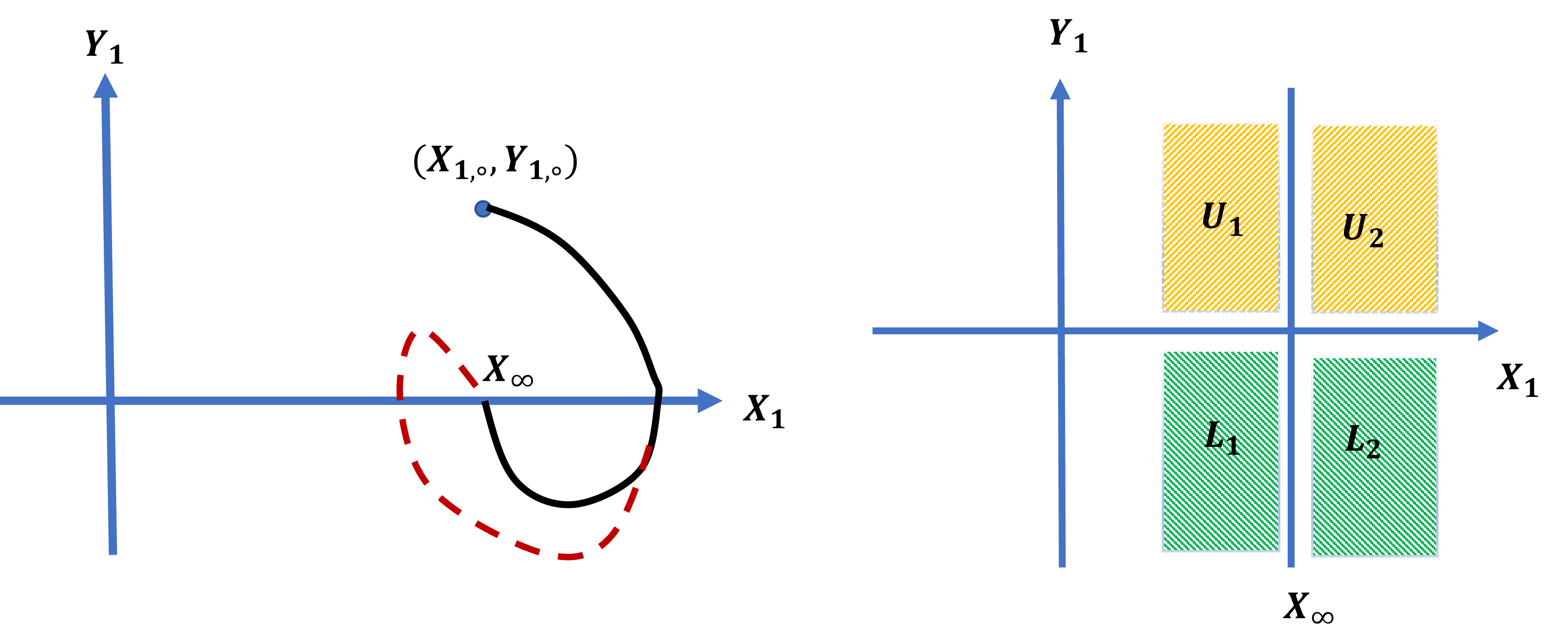}
    
        \caption{\small{The behavior of the trajectory for the different cases of initial values  $(X_{1, \circ}, Y_{1, \circ})$.}}
        \label{fig:9-10}
    \end{figure}
\end{proof}
Although in our setting $(X_{1, \circ}, Y_{1, \circ}) \in (0, X_\infty] \times \RR_+$, we need to understand the behavior of the dynamics when initial data is located in different positions around the equilibrium point. 
\begin{lemma}\label{T:negative_hit_zero}
    Suppose $(X_{1, \circ}, Y_{1, \circ}) \in L_1$. Then there exists a time $t_\circ \in (0, \MID_\infty)$ such that  $Y_1(t)>0$ for $t \in (t_\circ, \MID_\infty)$.     
\end{lemma}
\begin{proof}
    Since $Y_{1, \circ} <0$, the distance is decreasing. On the other hand, since the collision does not happen in the dynamic of $X_1$, there should be a time $ t_\circ \in (0, \MID_\infty)$ such that $Y_1(t_\circ) = 0$ and $\dot Y_1(t_\circ) >0$. Then from Lemma \ref{T:stay_positive}, $Y_1(t) >0$ for $t \in (t_\circ, \MID_\infty)$.
\end{proof}
As the result of the Lemma \ref{T:stay_positive}, either $Y_1(\MID_\infty) =0$ or we have that:
\begin{lemma}\label{T:U_2_case}
    Suppose, $(X_{1, \circ}, Y_{1, \circ}) \in U_2$. Then, there is a time $t_\circ < \MID_\infty$ such that $Y_1(t_\circ) = 0$ and $\dot Y_1(t_\circ) <0$. 
\end{lemma}
\begin{proof}
    Similar to previous lemmas, using the dynamics of $(X_1, Y_1)$, we have that $V(X_1) - \bar v >0$ since $X_{1, \circ} >X_\infty$ and therefore, $\dot Y_1 <0$. Hence the claim follows. 
\end{proof}
\begin{lemma}\label{T:final_lem}
    Suppose $(X_{1, \circ}, Y_{1, \circ}) \in L_2$. Then, $Y_1(t) <0$ for $t \in [0,\MID_\infty)$. 
\end{lemma}
\begin{proof}
    The proof is similar to the case of Lemma \ref{T:stay_positive} and by contradiction.  
    
    
\end{proof}
Therefore, if $(X_{1, \circ}, Y_{1, \circ}) \in L_2$, then either $Y_1(\MID_\infty) = 0$ or the case of Lemma \ref{T:negative_hit_zero} will be revisited. All of these cases will be repeated around the equilibrium point until the convergence happens.

    Employing the results of Lemma \ref{T:stay_positive}-\ref{T:final_lem}, and properties of the Hamiltonian \eqref{E:Hamiltonian}, we show that the rate of convergence of the trajectory $(X_1(t), Y_1(t))$ to the equilibrium point can be controlled by controlling the parameters $\alpha$ and $\beta$. 
    
    The following inequalities are the cornerstone of proving Theorem \ref{T:exp_stable_convergence}. 
\begin{lemma}\label{T:H_bound_expstab_condition}
        Let's define the domain $\mathcal C \Def (\delta_1 - X_\infty , \bar x - X_\infty) \times \RR$ (see \eqref{E:X1_bound} and \eqref{E:X1_LB} and the definition of $u,v$ before \eqref{E:new_hamiltonian}) which contains the equilibrium point. There exists constants $\underline k$ and $\overline k$ such that
        \begin{equation*}
            \underline k \norm{U}^2 \le H(u, v) \le \overline k \norm{U}^2
        \end{equation*}
        where $U = (u, v)^\sT \in \mathcal C$. 
    \end{lemma}
    \begin{proof}
        The right-hand side inequality is by Lipschitz continuity of function $V$ and the fact that $\bar v = V(X_\infty)$ and for $\overline k \Def \max \set{\tfrac 12, \alpha}$. To see the left-hand side of the inequality, we note that function $\tilde P$ (as in \eqref{E:new_hamiltonian}) is a convex function (see the illustration $P$ of Figure \ref{fig:geometry} and consider that the equilibrium point is shifted to the origin) and in addition, over the domain $\bar{\mathcal C}$, the closure, we have
        \begin{equation}\label{E:strong-convex-constant}
            \tilde P''(u) = \alpha V'(u + X_\infty) \ge \alpha k_{\eqref{E:strong-convex-constant}}  >0 
        \end{equation}
        for some constant $k_{\eqref{E:strong-convex-constant}}>0$. Therefore, $\tilde P$ is strongly convex on $\mathcal C$ which implies that 
        \begin{equation*}
            \tilde P(u) \ge \tilde P(0) + \tilde P'(0) u + \tfrac {\alpha k_{\eqref{E:strong-convex-constant}}}{2} u^2
        \end{equation*}
        Since $\tilde P(0) = \tilde P'(0) = 0$, we have that
        \begin{equation}\label{E:strong_convex}
            H(u, v) = \tfrac 12 v^2 + \tilde P(u) \ge \tfrac 12 v^2 + \tfrac {\alpha k_{\eqref{E:strong-convex-constant}}}{2} u^2 \ge \underline k \norm{U}^2
        \end{equation}
        for $\underline k  \Def \min \set{\tfrac 12, \tfrac {\alpha k_{\eqref{E:strong-convex-constant}}}{2}}$.
        This completes the proof.
    \end{proof}
    For the proof of Theorem \ref{T:exp_stable_convergence}, with a slight abuse of notation, we consider $\MID_\infty$ as in \eqref{E:equil_time} to denote the time that trajectory $t \mapsto u(t)$ approaches the origin (which is the equilibrium point here). 
    
    \begin{proof}[of Theorem \ref{T:exp_stable_convergence}]
        Let us fix $\eps >0$ such that $\ball((0,0)^\sT, \eps)$ be the region of attraction (for exponential stability) of the origin in the linearized model; see Remark \ref{R:global_asym_stab}. We recall Lemma \ref{T:stay_positive}-\ref{T:final_lem} and we define 
        \begin{equation}\label{E:Teps}
            \MID_\infty^\eps \Def \inf \set{t < \MID_\infty: u(t) < \eps/3}, 
        \end{equation}
       If for some values of $\alpha$ and $\beta$
        \begin{equation*}
            v(t) < \eps, \quad t \in [\MID_\infty^\eps, \MID_\infty),
        \end{equation*}
        i.e. is already in the domain of attraction, then the claim follows by exponential convergence of the linearized problem.
        
        Suppose on the contrary that for all values of $\alpha$ and $\beta$, $v(t) \ge \eps$ on $[\MID_\infty^\eps, \MID_\infty)$. Then over the domain $\mathcal C$
        \begin{equation*}
            \begin{split}
                \left(\alpha + \frac{\beta}{(u(t) + X_\infty)^2} \right)v^2(t) &\ge \tfrac 12 \left(\alpha + \frac{\beta}{\bar x^2}\right) v^2(t)\\
                & \qquad \qquad + \tfrac 12 \left(\alpha + \frac{\beta}{\bar x^2}\right) \eps^2 \\
                &\ge \tfrac 12 \left(\alpha + \frac{\beta}{\bar x^2} \right) \norm{(u(t), v(t))^\sT}^2
            \end{split}
        \end{equation*}
        for  $t \in [\MID_\infty^\eps, \MID_\infty)$, and 
        where the last inequality is by \eqref{E:Teps}. Using \eqref{E:new_hamiltonian}, we have that
        \begin{equation}\label{E:dotH_bound_asym}\begin{split}
            \dot H(u(t) , v(t)) & \le - K_{\eqref{E:dotH_bound_asym}}\norm{(u(t), v(t))^\sT}^2, \quad t \in [\MID_\infty^\eps, \MID_\infty).
        \end{split}\end{equation} 
        where, $ K_{\eqref{E:dotH_bound_asym}} \Def \tfrac 12 \left(\alpha +  \frac{\beta}{\bar x^2}\right)$.
   
    Using Lemma \ref{T:H_bound_expstab_condition} and \eqref{E:dotH_bound_asym}, we can write
    \begin{equation*}
        \dot H(u(t), v(t)) \le - \frac{K_{\eqref{E:dotH_bound_asym}}}{\overline k} H(u(t), v(t)),\quad t \in [\MID_\infty^\eps, \MID_\infty).
    \end{equation*}
    Using Gronwall's inequality, we get
    \begin{equation*}
        H(u(t), v(t)) \le H(u(\MID_\infty^\eps), v(\MID_\infty^\eps)) \exp\set{- \frac{K_{\eqref{E:dotH_bound_asym}}}{\overline k} t},
    \end{equation*}
    for $t \in [\MID_\infty^\eps, \MID_\infty)$. 
    Once more, using Lemma \ref{T:H_bound_expstab_condition}, we will have that 
    \begin{equation*}
        \begin{split}
            \norm{(u(t), v(t))^\sT}^2 & \le \frac{1}{\underline k} H(u(t), v(t)) \\
            & \le \frac{1}{\underline k} H(u(\MID_\infty^\eps), v(\MID_\infty^\eps)) \exp\set{- \frac{K_{\eqref{E:dotH_bound_asym}}}{\overline k} t}\\
            & \le \frac{ h_\circ}{\underline k}  \exp\set{- \frac{K_{\eqref{E:dotH_bound_asym}}}{\overline k} t}, \quad t \in [\MID_\infty^\eps, \MID_\infty),
        \end{split}
    \end{equation*}
    where the last inequality is from \eqref{E:energyLevel} and \eqref{E:hcircdef}. 
    But comparing $K_{\eqref{E:dotH_bound_asym}}$ and $\overline k$ shows that for sufficiently large values of $\alpha$ and $\beta$, $\norm{(u(t), v(t))^\sT} <\eps$ for some $ t\in (\MID_\infty^\eps, \MID_\infty)$ which contradicts our initial assumption. Therefore, the statement of the theorem follows
    (see Figure \ref{fig:alpha_beta_effect}).
    \begin{figure}
        \centering
        \includegraphics[width=0.9\columnwidth]{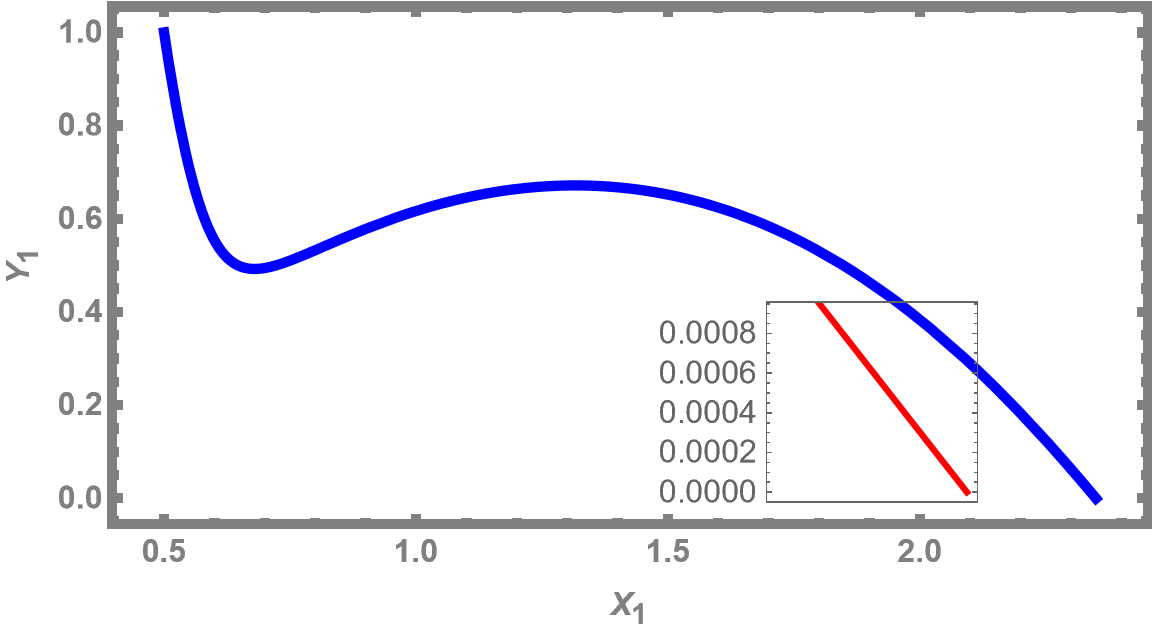}
        \includegraphics[width=0.90\columnwidth]{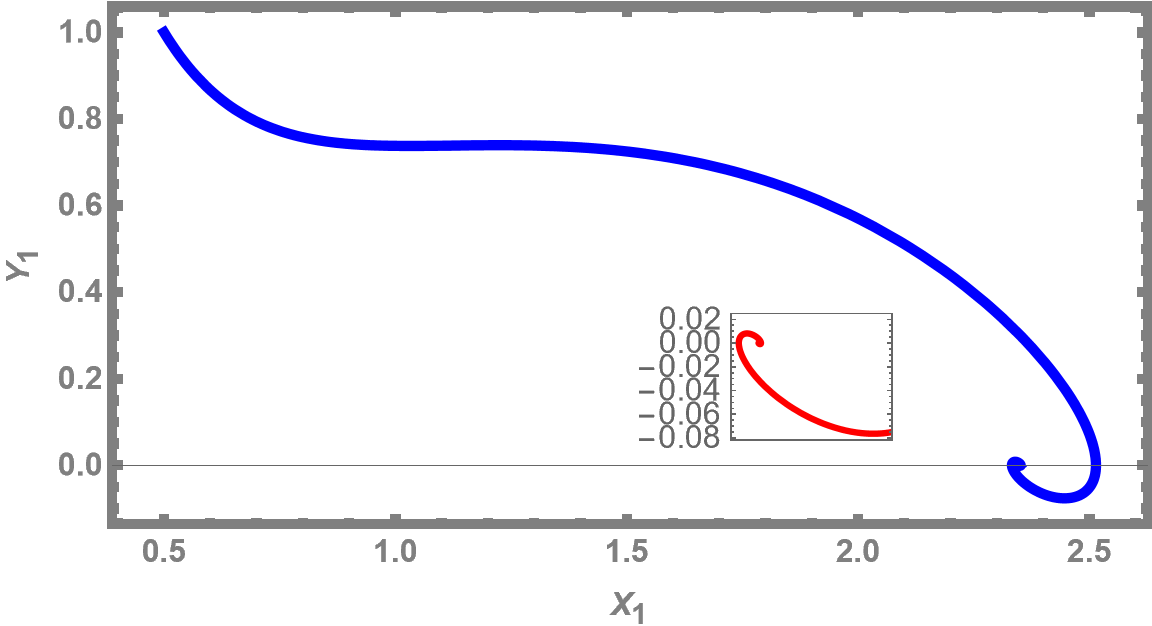}
        \caption{\small{The first figure is for $\alpha = 3$, $\beta = 2$ and it converges to the rest point. The second figure is with respect to $\alpha = \beta = 1$. Hence for relatively small $\alpha$ and $\beta$, the cases of Lemma \ref{T:U_2_case}-\ref{T:final_lem} happen. The red curves zoom into the behavior of the trajectories near the equilibrium point. Other parameters are $\bar v = 1.3$, $(X_{1, \circ}, Y_{1, \circ})=(0.5, 1)$.}}
        \label{fig:alpha_beta_effect}
    \end{figure}
     \end{proof}

\section{Dynamics of Other Vehicles}
The analysis of the properties of the dynamics of interaction between vehicle $n$ and $n-1$ for $n \ge 2$ requires an in-depth understanding of the interaction between vehicles $n-1$ and $n-2$ (the leading vehicles). In this section, using the results of section \ref{S:two_vechicles}, we will consider the interaction between vehicles $n$ and $n -1$ for $n= 2$ (the interaction between vehicles two and three). 

For the purpose of our analysis (in particular collision analysis), we need to work with the \textit{difference} flow $(X_2- X_1,Y_2 - Y_1)$ rather than the flow $(X_2, Y_2)$. In particular, $X_2 - X_1 = 0$ means collision. Therefore, it would be reasonable to introduce the change of variables
\begin{equation*}\begin{split}
    \xi_1 \Def X_1, \quad \xi_2 \Def X_2 - X_1 \\
    \zeta_1 \Def Y_1 \quad \zeta_2 \Def Y_2 - Y_1,
\end{split}\end{equation*}
and the difference dynamical model of \eqref{E:main_dynamics} then reads
\begin{equation}\label{E:differece_dynamics}
    \begin{cases}
        \dot \xi_1 = \zeta_1 \\
        \dot \zeta_1 = -\alpha \lb V(\xi_1) - \bar v \rb - \alpha  \zeta_1 - \beta \frac{\zeta_1}{(\xi_1)^2} \\
        \dot \xi_2 = \zeta_2 \\
        \dot \zeta_2 = - \alpha \lb V(\xi_2) - V(\xi_1) \rb -  \alpha \zeta_2 -  \beta \left(\frac{\zeta_2}{(\xi_2)^2} - \frac{\zeta_1}{(\xi_1)^2} \right)\\
        (\xi_1(0), \zeta_1(0)) = (\xi_{1, \circ}, \zeta_{1, \circ }) = (X_{1, \circ}, Y_{1, \circ})\\
        (\xi_2(0), \zeta_2(0))= (\xi_{2, \circ}, \zeta_{2, \circ}) = (X_{2, \circ} - X_{1, \circ }, Y_{2, \circ } - Y_{1, \circ }).
    \end{cases}
\end{equation}
First, we look at the existence of the solution of the dynamics of \eqref{E:differece_dynamics}. We define the state space 
\begin{equation}\label{E:new_domain}
\mathcal D \Def \left((0, \infty)\times \RR\right)^2 \ni (\xi_{1, 0}, \zeta_{1, 0}, \xi_{2, 0}, \zeta_{2, 0}) \end{equation}
and the flow $\Xi \Def (\xi_1, \zeta_1, \xi_2, \zeta_2)$. 
From the abstract theory of dynamical systems \cite{teschl2012ordinary}, the solution of \eqref{E:differece_dynamics} exists on a maximal interval $[0, \MID^c)$, for some $\MID^c >0$. 

It should be noted that if $\MID^c = \infty$ then the system is globally well-behaved and consequently no collision occurs. In what follows, we first show that if $\MID^c < \infty$, i.e. the solution does not exist at all times, $\MID^c$ should represent the time that collision happens in the system; i.e. $\xi_2(t) \to 0$ as $t \to \MID^c$. 

\textbf{Assumption.} We suppose that 
\begin{equation}\label{E:finite_collision_time}
\MID^c < \infty. \end{equation}
We start with the implications of \eqref{E:finite_collision_time}.
As $t \nearrow \MID^c$, the flow $\Xi$, either grows unbounded, or $\Xi \in \partial \mathcal D$. We recall that under the conditions of Theorem \ref{T:exp_stable_convergence}, $\zeta_1$ vanishes at $\MID_\infty$. Therefore, if $\MID^c \ge \MID_\infty$ then for $t \in [\MID_\infty, \MID^c)$ the dynamics $\dot \zeta_2$ in \eqref{E:differece_dynamics} will be the same as dynamics of $\dot \zeta_1$ and so the solution exists for all $t \ge \MID_\infty$; i.e. $\MID^c = \infty$. On the other hand, $\zeta_2 \nearrow \infty$ which implies $\xi_2 \nearrow \infty$ is prohibited since by properties of function $V$ and boundedness of $\xi_1$ and $\zeta_1$, this implies $\dot \zeta_2 \to -\infty$ which is a contradiction. Finally, $\zeta_2 \searrow -\infty$ implies that $\xi_2 \searrow - \infty$ by the third equation of the dynamical system \eqref{E:differece_dynamics} which is not admissible in the domain $\mathcal D$. Putting all together, for $\MID^c < \MID_\infty$, by Lemma \ref{T:stay_positive}-\ref{T:final_lem} and since $\xi_1(t) \notin \set{0}$ for $t < \MID_\infty$, from \eqref{E:new_domain} we must have 
\begin{equation}\label{E:collition_time_limit}
    \lim_{t \nearrow \MID^c} \xi_2(t) = 0,
\end{equation}
or in other words, if assumption \eqref{E:finite_collision_time} holds, then$\MID^c$ must be the collision time. 
Let's look at the implications of the finite collision time. Under such an assumption, there exists a time 
\begin{equation*}
    \check t \Def \sup \set{t < \MID^c: \xi_2 > \tfrac 12 \delta_2} ,
\end{equation*}
where we set
\begin{equation}
    \label{E:delta_2}
    \delta_2 \Def \min \set{\xi_{2, \circ }, \delta_1},
\end{equation}
where $\delta_1 < X_\infty$ is defined in \eqref{E:X1_LB}.
In other words, if the collision time is finite, then there should be a time $\check t$ after which the trajectory $\xi_2(t) \le \tfrac 12 \delta_2$, for $t \in [\check t, \MID^c)$.
We study the behavior of the $\dot \zeta_2$ in this region. The next result shows that in this region, $\zeta_2(t) <0$; i.e. the follower vehicle is moving faster than the leading one. 
\begin{lemma}\label{T:invar_set}
  The Set $\mathcal U \Def \set{(x, y): x \in (0, \tfrac 12 \delta_2), y \in \RR_+}$ is invariant with respect to the trajectory $t \in [\check t, \MID^c) \mapsto (\xi_2(t) , \zeta_2(t))$. In other words, if $\zeta_2(t) > 0 $ for some $t \in [\check t, \MID^c)$, then it must remain positive. \end{lemma}
\begin{proof}
Figure \ref{fig:region} illustrates the proof argument. Suppose on the contrary that $\zeta_2(t) <0$ for some $t \in [\check t, \MID^c)$. This implies, by definition of $\check t$, there exists a time $t_\circ$ such that $\zeta_2(t_\circ) =0$ and $\dot \zeta_2(t_\circ) <0$. But using the dynamics of $\zeta_2$ in \eqref{E:differece_dynamics} as well as \eqref{E:delta_2}, must have that 
\begin{equation*}
    \dot \zeta_2(t_\circ) = - \alpha \lb V(\xi_2(t_\circ)) - V(\xi_1(t_\circ) \rb + \beta \frac{\zeta_1}{(\xi_1)^2} >0, 
\end{equation*}
which is a contradiction. 
\begin{figure}
    \centering
    \includegraphics[width=2in]{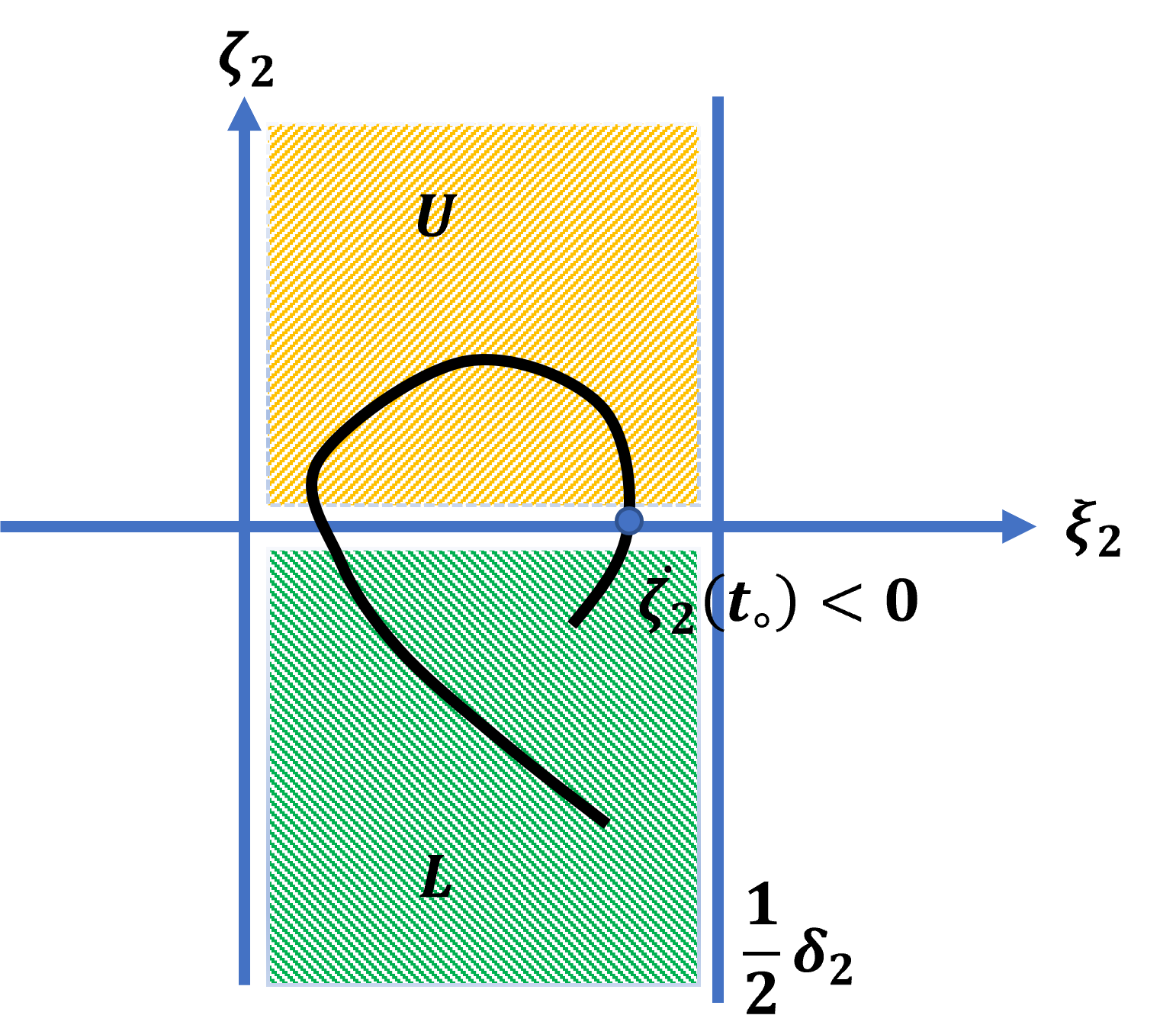}
    \caption{\small{Illustration of the proof discussions.}}
    \label{fig:region}
\end{figure}
\end{proof}
\begin{remark}\label{R:zeta_2_sign}
Lemma \ref{T:invar_set} along with the assumption \eqref{E:finite_collision_time} shows in particular that $\zeta_2(t) <0$ for $t \in [\check t, \MID^c)$. \end{remark}

\begin{figure}
    \centering
    \includegraphics[width=2.5in]{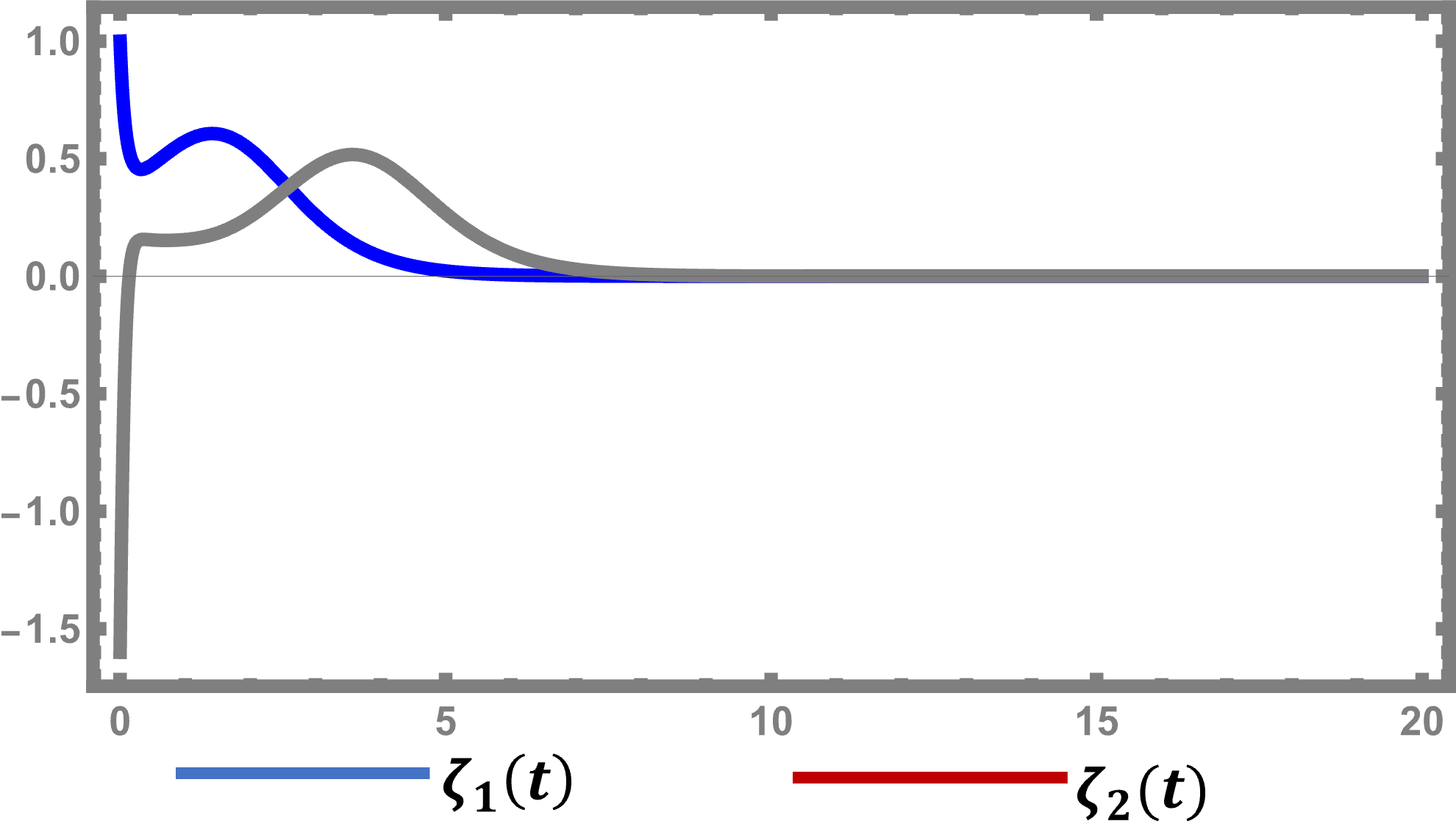}
    \includegraphics[width=2.35in]{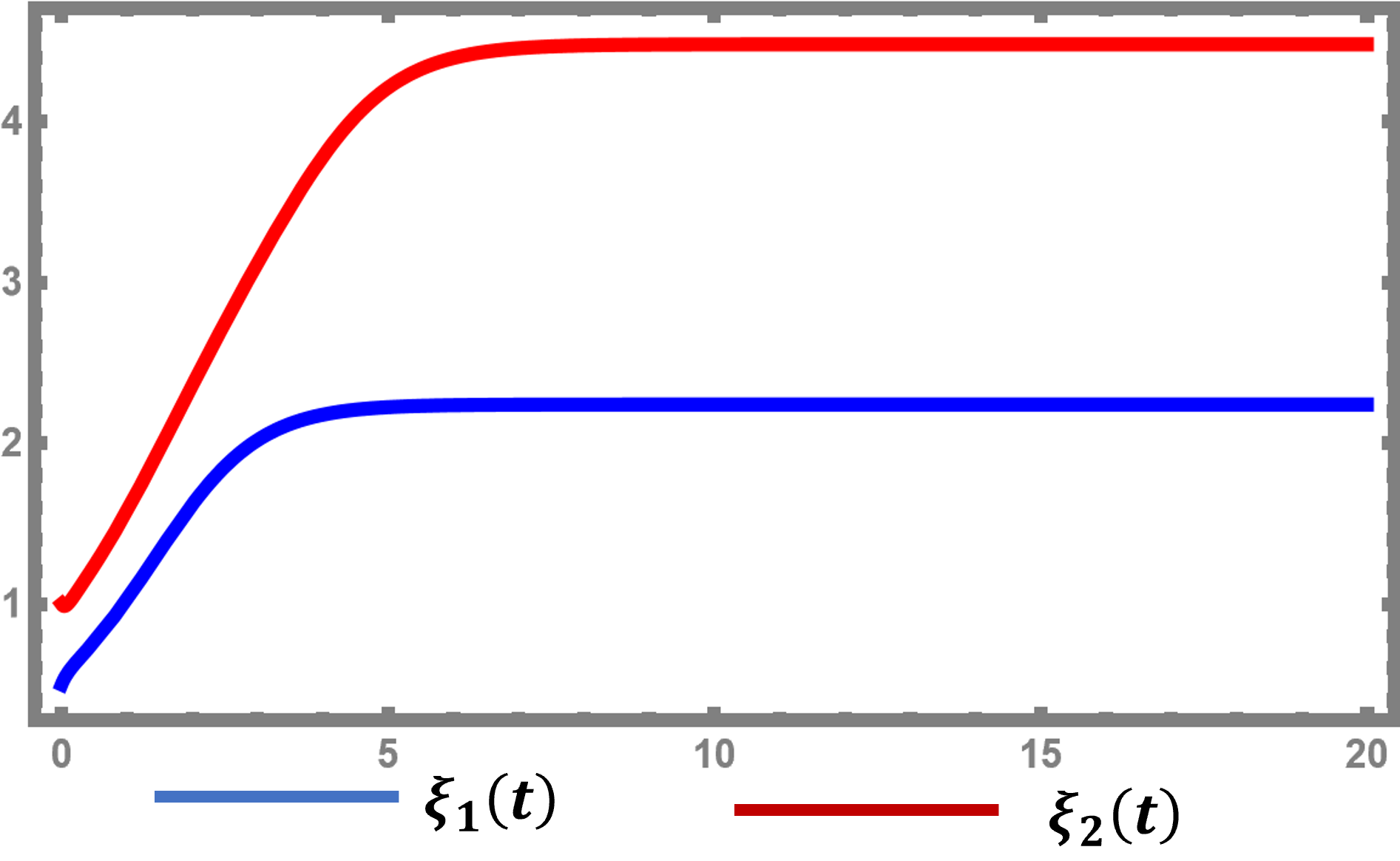}
    \includegraphics[width=0.83\columnwidth]{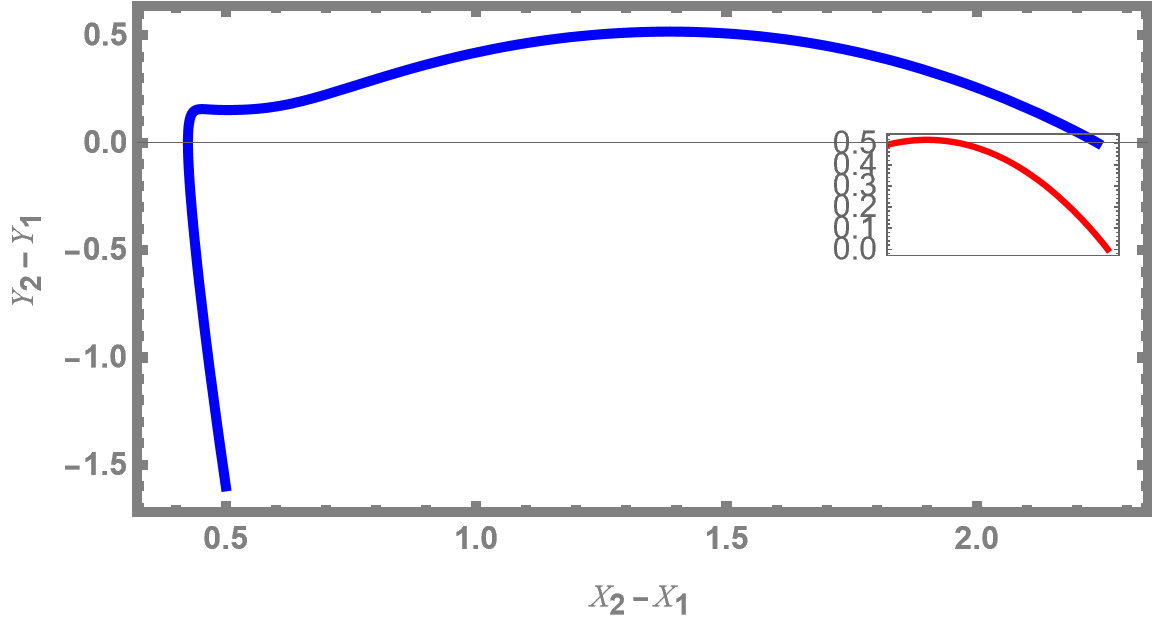}
    \caption{\small{Illustration of the trajectory and orbit of the dynamics between the first two and the second two vehicles. The red curve in the last illustration shows the behavior of the orbit near the rest point.}}
    \label{fig:18-20}
\end{figure}

Next, we show that while on $[\check t, \MID^c)$, $\xi_2$ is small and $\zeta_2 <0$, deceleration will be encouraged, which is the mere hope to prevent the collision. In other words, if deceleration is not strong enough, the collision will happen.  
\begin{lemma}\label{T:positive_accel}
 For $t \in [\check t, \MID^c)$, we have that 
    \begin{equation*}
        \dot \zeta_2(t) > 0.
    \end{equation*}
\end{lemma}
\begin{proof}
    Let's consider the dynamics of $\dot \zeta_2$ in \eqref{E:differece_dynamics}. Then, \eqref{E:delta_2} implies that $- \alpha \lb V(\xi_2) - V(\xi_1) \rb >0$. Remark \ref{R:zeta_2_sign} show that $- \alpha \zeta_2 >0$. Finally $\zeta_1 >0$ implies that the last term should also be positive and this completes the proof. 
\end{proof}
Now, we are ready to show that our assumption \eqref{E:finite_collision_time} and its implications (the previous results) lead to a contradiction, i.e. the collision cannot happen in finite time. 
\begin{proposition}\label{T:collison}
Under the conditions of Theorem \ref{T:exp_stable_convergence}, $\mathcal T^c = \infty$. In other words, collision does not happen in the dynamical model \eqref{E:differece_dynamics}.
\end{proposition}
\begin{proof}
Using \eqref{E:collition_time_limit}, the idea of the proof is to show that $\xi_2(t) > 0$ for $t \in [\check t, \MID^c)$ which implies the statement of the theorem. The result of Lemma \ref{T:positive_accel} (monotonicity of $\zeta_2(t)$) suggests that, given the trajectory of $t \mapsto (\xi_1(t), \zeta_1(t))$, we should be able to locally write
\begin{equation}\label{E:Psi_dynamics}
    \xi_2(t) = \psi(\zeta_2(t)), \quad \text{for $t \in [\check t, \mathcal T^c)$},  
\end{equation}
for some function $\psi$ which will be constructed below. Employing \eqref{E:differece_dynamics}, we have that
\begin{equation}\label{E:transit}
\begin{split}
    \zeta_2 = \dot \xi_2 &= \psi'(\zeta_2) \dot \zeta_2 \\
        & = \psi'(\zeta_2) \Bigg \{ - \alpha \lb V(\psi(\zeta_2)) - V(\xi_1) \rb - \alpha \zeta_2 \\
        &\qquad \qquad \qquad    - \beta \left(\frac{\zeta_2}{(\psi(\zeta_2))^2} - \frac{\zeta_1}{(\xi_1)^2} \right) \Bigg\}.
\end{split}\end{equation}
Furthermore, thanks to strictly monotone behavior, the function $\zeta_2:[\check t, \mathcal T^c) \to [\check \zeta, \hat \zeta)$, where $ \check \zeta \Def \zeta_2(\check t)$ and $\hat \zeta \Def \zeta_2(\mathcal T^c)$, is a diffeomorphism. Let $\theta \Def \zeta_2^{-1}$, the inverse function of $\zeta_2$. Then, function $\zeta_1$ on $[\check t, \mathcal T^c )$ can be presented as the smooth function $\zeta_1 \circ \theta$ on $[\check \zeta, \hat \zeta)$ if $\mathcal T^c < \mathcal T_\infty$ and zero otherwise. A similar argument holds true for $\xi_1 \circ \theta$. \\
Let us now formalize the construction of $\psi$ by extending the function $\theta$ smoothly on the domain $(- \infty, 0)$ and defining a function
\begin{equation*}
{\scriptstyle
   \mathbf g(\zeta, \psi) \Def 
    \frac{-\zeta}{\alpha \lb V(\psi) - V(\xi_1(\theta(\zeta))) \rb + \alpha \zeta + \beta \left(\frac{\zeta}{\psi^2} - \frac{\zeta_1(\theta(\zeta))}{(\xi_1(\theta(\zeta)))^2} \right)}}
\end{equation*}
for $(\zeta, \psi) \in (- \infty, 0) \times (0 , \delta_2)$, and $(\zeta_1, \xi_1) \in (0, \bar y) \times(\delta_1, X_\infty)$. 
Therefore, using \eqref{E:transit}, the dynamical model can be presented by
\begin{equation}\label{E:phi_psi_dynamics}
\begin{cases}
   \psi'(\zeta)= \mathbf g(\zeta, \psi(\zeta)) \\
    \psi(\Check \zeta) = \xi_2(\check t)
\end{cases}\end{equation}
where $\check \zeta \Def\zeta_2(\check t)$ and $\xi_2(\check t) = \tfrac 12 \delta_2$. Through such construction, the dynamics of \eqref{E:phi_psi_dynamics} is well-defined and has a maximal interval of existence $(\mu_-, \mu_+) \subset (-\infty, 0)$ and contains the initial value $\check \zeta$. The construction \eqref{E:phi_psi_dynamics} creates a barrier dynamics through comparison with which we can show $\xi_2 >0$ on $[\check t, \MID^c)$ (see \eqref{E:Psi_dynamics}). 
\begin{theorem}\label{T:collision_avoid}
We have that 
\begin{equation*}
    \inf_{\bar \zeta \in {[\check \zeta, \mu_+)}} \psi(\bar \zeta) >0.
\end{equation*}
\end{theorem}
\begin{proof}
Let's consider the definition of $\psi'(\zeta)$ in \eqref{E:phi_psi_dynamics}. We recall that $\zeta_1 >0$ for $t \in [\check t, \MID^c)$, and by construction $\psi(\zeta) <  \delta_2 \le \delta_1 < \xi_1$. This implies that 
\begin{equation*}
    \psi'(\zeta) < 0, \quad \zeta \in [\check \zeta, \mu_+).
\end{equation*}
Dividing both sides of \eqref{E:phi_psi_dynamics}, we will have 
\begin{equation*}
    \begin{split}
        & \frac{\psi'(\zeta)}{\psi^2(\zeta)} = \\
        &   (- \zeta) \Bigg\{\alpha \psi^2(\zeta) \lb V(\psi(\zeta)) - V(\xi_1(\theta(\zeta))) \rb \\
        & \quad  + \alpha \zeta \psi^2(\zeta) + \beta \left(\zeta - \zeta_1(\theta(\zeta)) \frac{\psi^2(\zeta)}{(\xi_1(\theta(\zeta)))^2} \right) \Bigg\}^{-1}
    \end{split}
\end{equation*}
The right-hand side of this equation is in fact
\begin{equation}\label{E:psi_bound}
    \text{RHS} = \frac{- \zeta}{\text{Negative Terms + $\beta \zeta$}} > \frac{- \zeta}{\beta \zeta} = - \frac 1 \beta
\end{equation}
on $[\check \zeta, \mu_+)$. In addition, we define
\begin{equation*}
    \mathcal R(u) \Def \lb (\tfrac 12 \delta_2)^{-1} + \frac 1\beta (u - \check \zeta) \rb ^{-1}, \quad u \in [\check \zeta, 0).
\end{equation*}
Then, 
\begin{equation}\label{E:mathcal R_relations}
    \mathcal R(\check \zeta) = \tfrac 12 \delta_2 = \psi(\check \zeta), \quad \frac{\dot{\mathcal R}(u)}{\mathcal R^2(u)} = -\frac 1\beta , \quad u \in [\check \zeta, 0). 
\end{equation}
To compare $\mathcal R(\cdot)$ and $\psi(\cdot)$ which will help us showing the lower bound for $\psi$, we define
\begin{equation*}
    \begin{split}
        \Lambda (u) &\Def \left(\frac{1}{\psi(u)} - \frac{1}{\mathcal R(u)}\right)^+ \\
        & = \left(\frac{1}{\psi(u)} - \frac{1}{\mathcal R(u)}\right) \bOne_ {\lb \nicefrac{1}{\psi(u)} > \nicefrac{1}{\mathcal R(u)} \rb}\\
    \end{split}
\end{equation*}
for $u \in [\check \zeta, \mu_+)$. Therefore, if $\psi(u) < \mathcal R(u)$, we get
\begin{equation}\label{E:Lambda_negative_growth}
    \dot \Lambda(u) = \frac{\dot{\mathcal R}(u)}{\mathcal R^2(u)} - \frac{\psi'(u)}{\psi^2(u)} = - \frac 1 \beta - \frac{\psi'(u)}{\psi^2(u)} \le 0,
\end{equation}
and the last inequality is by \eqref{E:psi_bound}. But then from \eqref{E:mathcal R_relations} we have that $\Lambda(\check \zeta) = 0$, and hence \eqref{E:Lambda_negative_growth} implies that $\Lambda(u) \le 0$ for $u \in [\check \zeta, \mu_+)$ which is a contradition with the assumption $\psi(u) < \mathcal R(u)$. Therefore, 
\begin{equation}
    \inf_{u \in [\check \zeta, \mu_+)} \psi(u) \ge \inf_{u \in [\check \zeta, \mu_+)} \mathcal R(u) \ge \lb (\tfrac 12 \delta_2)^{-1} + \frac {\check \zeta}{\beta} \rb ^{-1}.
\end{equation}
This completes the proof. 
\end{proof}

We conclude the proof of proposition \ref{T:collison}, by showing that $\mu_+ = \hat \zeta = \zeta_2(\mathcal T^c)$. In other words, we show that the result of the Theorem \ref{T:collision_avoid} holds true for all $[\check t, \mathcal T^c)$ which in turn implies that $\mathcal T^c = \infty$. To do so, as mentioned before, $\zeta_2: [\check t, \mathcal T^c) \to [\check \zeta, \hat \zeta)$ is a homeomorphism. Therefore, we should have 
\begin{equation*}
    \theta \bigmid_{[\check t, \mathcal T^c)} \left([\check \zeta , \mu_+) \cap [\check \zeta, \hat \zeta) \right) = [\check t, \mathcal T'),
\end{equation*}
for some $\mathcal T' \le \mathcal T^c$ and we recall that $\theta = \zeta_2^{-1}$. From \eqref{E:phi_psi_dynamics} $(\psi(\zeta), \zeta)$ solves the dynamics of \eqref{E:differece_dynamics} and hence must coincide with $(\xi_2, \zeta_2)$ on $[\check t, \mathcal T')$. Let's assume that $\mathcal T' < \mathcal T^c$. In this case, by the extensibility theorem of dynamical systems, $(\xi_2(\mathcal T'-), \zeta_2(\mathcal T'-))$ should hit the boundary of $(0, \tfrac 12 \delta_2) \times (-\infty, 0)$ which is precluded by the fact that $\zeta_2 <0$ and $\dot \zeta_2 >0$ on $[\check t, \MID')$. Therefore, $\mathcal T' = \mathcal T$. However, in this case, we must have that 
\begin{equation*}
    \lim_{t \nearrow \mathcal T^c} \psi(\zeta) = \lim_{t \nearrow \mathcal T^c} \xi_2 = 0, 
\end{equation*}
which is prohibited by Theorem \ref{T:collision_avoid}. This contradicts our main assumption \eqref{E:finite_collision_time}. 
\end{proof}
To summarize, Proposition \ref{T:collison} states that the collision does not happen. Figure \ref{fig:18-20} illustrates the interaction of three vehicles in the form of the difference dynamics.

\section{Conclusions and future works}
In this paper, we presented a rigorous boundary layer analysis of the OVFL dynamical model near collision. Such analysis provides an in-depth understanding of the behavior of the dynamics (e.g. in a platoon of connected autonomous vehicles) especially when the system is forced out of equilibrium (for instance as a result of some perturbation in the system). Understanding the interaction of the singularity and behavior of the dynamics near collision is fundamental both from a theoretical standpoint and in designing efficient systems, such as adaptive cruise controls.  

This paper can be extended on several fronts. The theory can benefit from a broader definition of Hamiltonian which serves as a Lyapunov-type function to explain the boundedness and stability of the equilibrium solution. Utilizing this, further analysis is required to generalize the results in a rigorous way.

\bibliographystyle{IEEEtran}
\bibliography{reference}
\end{document}